\documentclass[11pt]{amsart}

\usepackage{graphicx} % [demo] option not to require actual image files 
\usepackage{caption}

\usepackage{subcaption}

\input xypic \xyoption{curve}

\newcommand{\B}[1]{{\mathbf#1}} % bold math
\newcommand{\C}[1]{{\mathcal#1}} % Calligraphic
\theoremstyle{plain}

\theoremstyle{definition}

\theoremstyle{definition}

\theoremstyle{remark}
\newtheorem{rem}{Remark}[section]

\usepackage{fancyhdr}
\begin{document}

% ****************************************************************************
% 				Title
% ***************************************************************************
\title{A Note on the Statistics of Riemann Zeros}         
\author{Lucian M. Ionescu}
\address{Department of Mathematics, Illinois State University, IL 61790-4520}
\email{lmiones@ilstu.edu}
\date{\today}         

%==================================================
%                      		Abstract
%==================================================
\begin{abstract}
Evidence of an algebraic/analytic structure of the Riemann Spectrum,
consisting of the imaginary parts of the corresponding zeros, is reviewed,
with emphasis on the distribution of the image of the primes
under the Cramer characters $X_p(t)=p^{it}$.

The duality between primes and Riemann zeros,
expressed traditionally as the Riemann-Mangoldt exact equation, 
is further used to investigate from a statistical point of view, the correspondence 
between the POSet structure of prime numbers and this yet unknown 
structure of R-Spec.

Specifically, the statistical correlation coefficient $c(p,q)=<X_p,X_q>$
is computed, noting ``resonances'' at the generators $q$ of the 
symmetry group $Aut_{Ab}(F_p)$ of finite field $F_p$.

A program for further studying the Riemann zeros from a pro-algebraic point of view,
is presented.
\end{abstract}

\maketitle
\setcounter{tocdepth}{3}
\tableofcontents

%***************************************************************************************
% 			Introduction
%***************************************************************************************
\section{Introduction} 
The Riemann Hypothesis (RH) is one of the most important problem of modern mathematics
not only historically speaking \cite{MilleniumP}, but also because 
it is a ``grand junction'' of several areas in Mathematics and Physics \cite{LI:Rem}, 
and with impact on thousands of conditional theorems.
To prove it requires, it seams, a better understanding of what the Riemann zeros $\rho_n$ are,
with respect to their duality to prime numbers \cite{Mazur-Primes}, 
expressed traditionally as the Riemann-Mangold exact formula \cite{Garrett-REF},
which can be interpreted as a multiplicative Poisson-like Trace Formula,
or arguably better, as a Fourier Transform in the sense of distributions \cite{LI:PNRZ}.

What is in fact to be better understood is the {\em Riemann Spectrum} 
$R-Spec=\{t_n|\rho_n=1/2 +it\}$, i.e. the imaginary parts of the zeros, 
supposed to be real, according to the RH; since computationally
we have large samples of such zeros, with as needed accuracy, 
and we can prove by counting arguments they satisfy the RH.

In this article we focus on what we will call {\em Cramer characters} $s^iz$
pairing $R-Spec$ and $Z-Spec$, the rational primes.

Several studies, e.g. by Rademacher \cite{Rademacher}, Ford and Zaharescu \cite{Ford1},
were concern with the distribution of these ``matrix coefficients'' $X_p(t)=p^{it}$ (periods? \cite{KZ}),
under a different guise: truncation to the unit interval of $log(p)t$ \cite{Ford1}.
That these are the ``players'' in the Primes-Zeros duality \cite{Mazur-Primes}, 
is clear: $p^{it}=e^{2\pi i (\log(p)t/2\pi)}$ \cite{LI:PNRZ}.

On the other hand, a relatively new structure of the rational primes was emphasized 
in \cite{LI:POSet}, previously ``confined'' to the primality tests work by Pratt \cite{Pratt} 
(prime certificates and what the author now calls {\em Pratt trees}). 
This relation between prime numbers reflects the symmetries of finite fields,
a kind of resolution from objects to generators and relations, 
fundamental in mathematics, with deep historical roots in cohomology theory.

It is natural therefore to compare the two: 
the know, yet not much studied, POSet structure (the $Aut_{Ab}(F_p)$ ``loop functor''),
and the unknown R-Spec structure, via the not so well documented at the level of proofs
\footnote{A detailed proof of Fourier duality at the level of distributions
is not known to the author; see also \cite{Marco}.},
Prime-Zeros duality.

\cite{Ford1} studies the distribution of $X_p$,
which we will interpret in a statistic flavor as a {\em random variable}
(see also \cite{LI:PNRZ}, \S2, p.5).
In a sequel \cite{Ford2}, a 2-dimensional distribution was studied, 
in order to see possible $(X_p,X_p')$ correlations.

In this article we investigate the {\em correlation coefficient} between $X_p$ and $X_q$,
where $q|p-1$, i.e. searching for a ``hierarchic correlation'', 
complementing their search for a direct algebraic dependence between Riemann zeros.
A resonance phenomenon is apparent, yet no algebraic theoretic explanation is yet available
for the author.

Recall that there is evidence for such a structure on $R-Spec$:
the ``Graham points entropy'' is low \cite{Entropy}, and 
``L-function zeros know of each others'' \cite{Marco}.
It is as if there is an algebraic basis (crystal basis? See {\em Quasicrystals}, \cite{LI:PNRZ}),
and the zeros are cofinal \cite{Marco}: any infinite subset determines the hole set.

In conclusion, a further combined investigation is required, 
computational and theoretic,
including an adelic interpretation of the duality,
from the point of view of algebraic quantum groups \cite{VanDaele}. 
Since quantum crystal bases seam to be a topic to ponder upon in these context,
categorifying the duality towards a Tannaka-Krein duality of the category of 
finite abelian groups explaining the above prime-zeros duality is a direction to be considered.

%***************************************************************************************
% 			Primes and zeros
%***************************************************************************************
\section{The Primes and Zeros Duality}
To prove the Prime Number Theorem regarding the distribution of primes, 
led Riemann to study the Riemann Zeta Function, relating the primes and its non-trivial zeros
in the Riemann-Mangoldt exact formula \cite{Mazur-Primes,Garrett-REF}.
This is a global indication of a duality between primes and zeros,
allowing to express {\em globally}, as an ensemble,
the primes in terms of zeros (Fig.1):
\begin{figure}
\includegraphics[width=3in, height=2in]{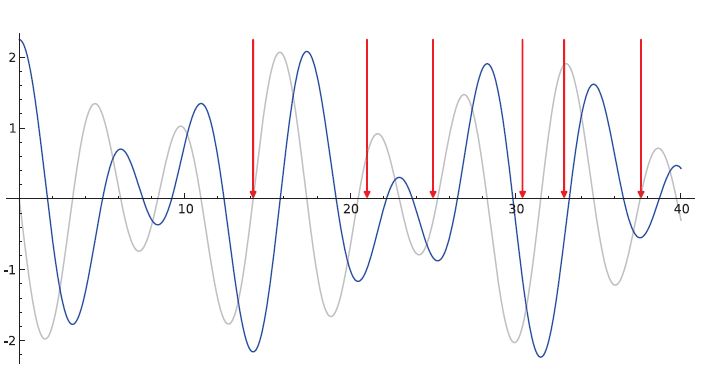}
\caption{From \cite{Mazur-Primes}, part IV, p.107}
\label{Mazur-diagram1}
\end{figure}
where the peaks seem to converge to the zeta zeros $\theta_n$,
and conversely, the zeros in terms of prime numbers (Fig.2):
\begin{figure}
\includegraphics[width=3in, height=2in]{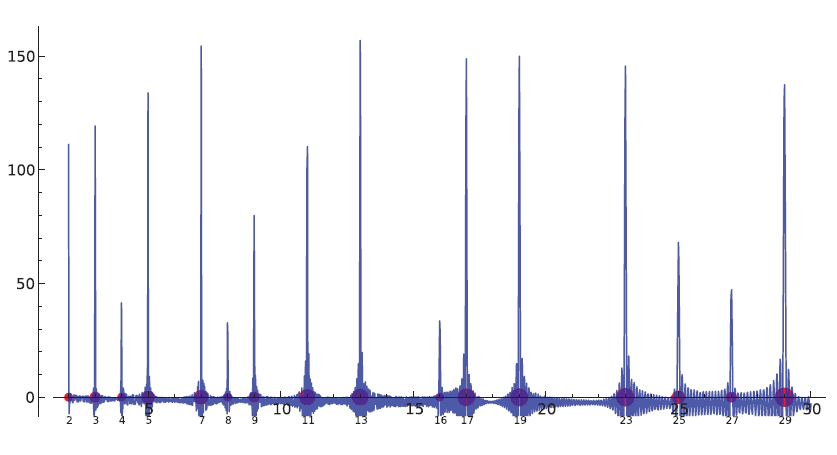}
\caption{From \cite{Mazur-Primes}, part IV, p.117}
\label{Mazur-diagram2}
\end{figure}

The corresponding sums  (see \cite{LI:PNRZ} for additional details),
may further be interpreted as Fourier Transform (FT) in the sense of distributions:
\begin{equation} \label{E:FT-DiracCombs}
FT(D_{Z-Spec)}=\C{D}_{R-Spec},
\end{equation}
i.e. the {\em Dirac combs} of the two spectra are Fourier dual
\footnote{A well documented proof is still missing, it seams; see also \cite{Marco}.}.

Since  {\em new structures} were uncovered on both spectra, namely
the {\em POSet structure on primes} \cite{LI:POSet} and 
non-trivial distributions of $R-Spec$ \cite{Ford1, Ford2},
it is natural to attempt to relate the two.
This article is mainly concerned in what follows with a 
statistical study of $R-Spec$, 
with hindsight from the POSet structure on the ``other side'' of the duality.

%***************************************************************************************
% 			Statistics of R-Spec
%***************************************************************************************
\section{Statistical Exploration of R-Spec}
Adopting a statistical interpretation of our lack of information regarding the 
yet unknown (pro)algebraic structure of the Riemann zeros, 
hinted for instance by its low Graham points entropy \cite{Entropy},
we will treat the $R-Spec$ as our {\em population},
with a typical {\em sample} of say the first few $100,000$ zeros, as our population
\footnote{The RH is of course assumed, so that Riemann zeros $\rho=1/2+it$ correspond to 
the imaginary parts $t$, forming the Riemann spectrum $R-Spec\subset \B{R}$, as denoted in the introduction.}.
The formula counting the number of imaginary parts (zeros) in the imaginary axis:
$$N(T)=\frac{T}{2\pi}(\log \frac{T}{2\pi}-1)+O(log T),$$
is ``reciprocal'' to the corresponding counting function of primes ({\em Prime Number Theorem} \cite{Mazur-Primes}.
\begin{rem}
This formula has a {\em quantity of information} flavor when compared asymptotically to:
$$I(y)=\int_0^y log(x)dx, \quad y=T/2\pi.$$
\end{rem}
In this way $(R-Spec, d\mu)$ is our primary (non-compact) measure space. 
% The RV's
Each prime number $p\in Z-Spec$ determines a {\em random variable} $X_p$,
also referred to as the p-sector of the Riemann Spectrum:
$$X_p:R_{Spec}\to C^\times, \quad X_p(t)=p^{it}.$$
Its values belong to the unit circle, in the complex plane.
\begin{rem}
Algebraically, they are entries of the matrix of the bi-characters 
$s^{iz}: (C^\times,\cdot) \times (C,+)\to (C^\times,\cdot))$,
in the corresponding bases $R-Spec$ and $Z-Spec$
% Long footnote!
\footnote{$Z-Spec$ is a basis of $(Q,\cdot)$, but $R-Spec$ is yet unclear 
in what sense it is such a basis; possibly in the context of 
the Fourier duality of the rationals as an {\em algebraic quantum group} \cite{VanDaele,LI:QQG}.}
\end{rem}

Exponentiation $p^{it}$ of the imaginary part $t$ of the Riemann zero $\rho=1/2+it$
corresponds to the truncation of the imaginary part to the unit interval, after rescaling by $2\pi$,
as in the work by Ford and Zaharescu \cite{Ford1,Ford2}.

% ***************   Distributions of p-Sectors X_p  **************************
\subsection{The distribution of p-sector $X_p$}
While the distribution of the $R-Spec$ after reduction to unit interval after scaling by $(\log p)\2\pi$, 
or to the unit circle via exponentiation, is uniform,
the ``$p$-sector measurements'' via the RV $X_p$ exhibits non-trivial distributions. 
For example, see Fig. 3, where $p=2$, and Fig. 4, with $p=5$.
% Replace on the same line, maybe?
\begin{figure}[h!]
\includegraphics[width=2in,height=1in]{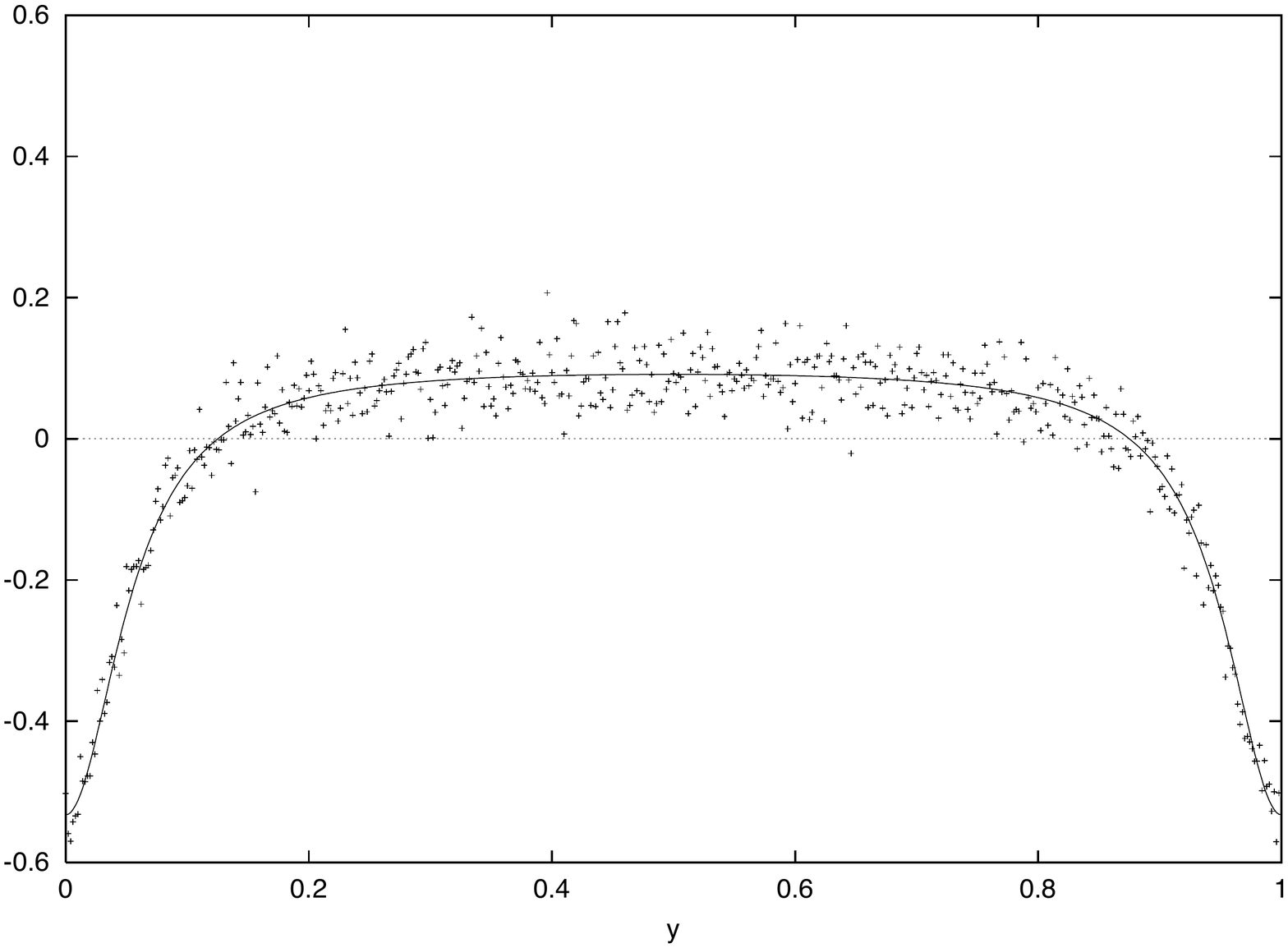} % {Fig1--2.png}
\caption{$p=2, \ q=1$}
\end{figure}
Similar non-trivial distributions occur for $p=5$ and $q=3$ (interpreted as a compression factor),
shown in Fig.4.
\begin{figure}[h!]
\includegraphics[width=4in,height=2in]{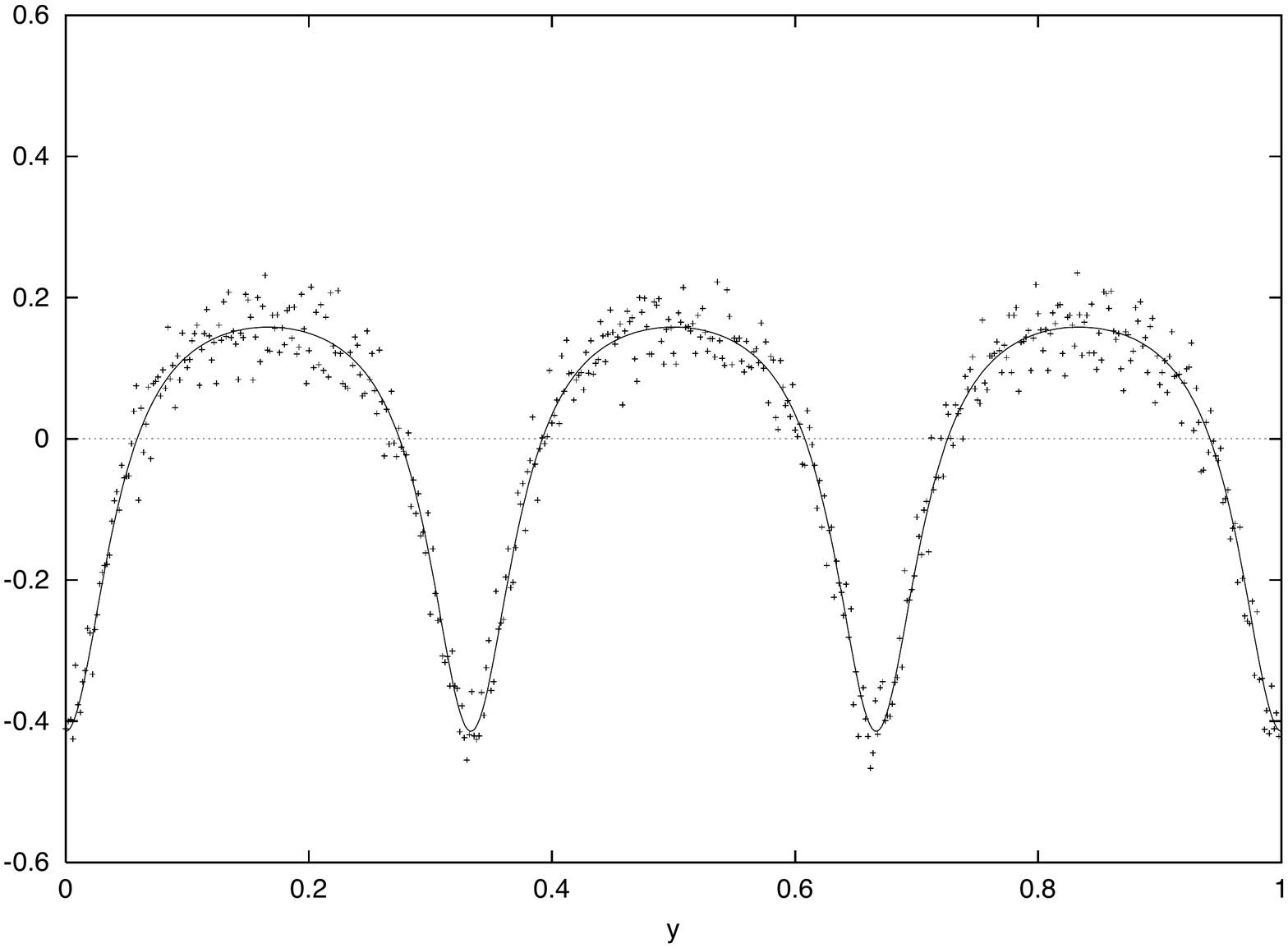} % {Fig1--5.png}
\caption{$p=5, \ q=3$}
\end{figure}
The role of the compression factor $q=3$ present in Fig. 4, 
leading to ``harmonics'', together with the definition of the frequency of zeros 
as well as the definition of the mathematical model of the density,
can be found in \cite{Ford1}.

\subsection{The bi-distributions $X_{p.q}$}
A further study of multivariable distributions 
corresponding to vectors $(p_1, ...,p_n)$ of such scaling factors 
$\log P=(\log p_i, ...,\log p_n)$ was reported in \cite{Ford2}.
Interesting ``lattice-like'' diagrams were obtained,
for example with $p_1=2, p_2=3$ as in Fig. 5,
after a linear integral transformation $M\alpha=\log P/2\pi$,
with $M=\left[\begin{matrix}1 & 1\\ 1 & -1\end{matrix}\right]$ 
(see \cite{Ford2} for details).
\begin{figure}[h!]
\includegraphics[width=4in,height=2in]{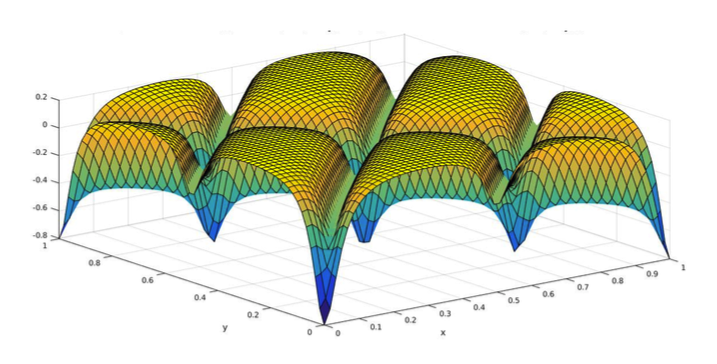}
\caption{$p_1=2$, $p_2=3$}
\end{figure}
\begin{rem}
A few questions regarding the above study are, perhaps, in order.

The distributions $g_\alpha(X)$ Eq. (1.3), \cite{Ford2}, p.2, 
seams to correspond to the primes-zeros Fourier sums (duality): 
\cite{Mazur-Primes}, Ch.30; \cite{LI:PNRZ}, \S3.2;
are these a finite dimensional truncation of the global Fourier Transform,
converging to Dirac delta functions?

The integral transformations of the lattice of primes (Dirac comb of $\log p$),
may in particular belong to the modular group.
In the context of Poisson trace formula \ref{E:FT-DiracCombs},
is there some kind of invariance involved?

Finally, when fractional exponents of primes are involved, as in Example 1 or 2 (loc. cit.),
should an algebraic study be done in the context of number fields,
adjoining roots of $x^n=p$? (Ramification Theory).
\end{rem}
In what follows we will study the correlation between two such RV $X_p$ and $X_q$,
and note a remarkable ``resonance'', i.e. high correlation coefficient, 
when $p$ and $q$ are ``related'' in the POSet sense: when $q|p-1$,
i.e. the prime $q$ corresponds to a {\em symmetry} of the 
{\em discrete Klein-Galois Geometry} of the finite field.
% Discrete Vector Spaces
\begin{rem}
Such a ``discrete space of vectors'' (because we can add and scale elements),
that we normally call abelian group or $Z$-module,
when ``localized, becomes $Q-mod$, i.e. a (rational) vector space).
\end{rem}

\subsection{The POSet Correlation}
Primes do not come in isolation, either \cite{LI:POSet}.
\subsubsection{The Complexity of Primes}
% Brief
Primes (and their powers) are but ``shadows'' of finite elementary abelian groups.
As discrete ``spaces of vectors'', they have symmetries in the sense of Klein geometry. 
The larger the factorization of $p-1$ the more symmetries the corresponding geometry has.
Sophie Germain primes $p=2q+1$ are especially simple, 
explaining her success in handling the reduction
of Fermat problem to such a ``simpler'' finite characteristic.
To build more complex primes $p>>q$, in the sense of the POSet of primes,
with $p$ ``on top'' of $q$, i.e. having (multiplication by) $q$ as a symmetry,
use {\em Euclid's trick} $p=2q_1...q_r+1$, and of course a ``quality control'' primality test.

Now is there a ``mirror reflection'' of this relation / hierarchy of primes' complexity,
at the level of the Riemann Spectrum!
There is only one way to find out, fast that is: computing the statistical correlation coefficient, next,
and if there is enough evidence, look for a conjecture etc. 

\subsubsection{Relations with $R-Spec$}
A preliminary theoretical study of the relation between primes and zeros,
from an adelic point of view, is contained in \cite{LI:PNRZ}.
This article is a computational follow-up study,
of the statistics correlation of two p-sectors of the Riemann spectrum
$c(X_p,X_q)$, when $<<p$ are related in the POSet of primes.

The programs were written in SAGE/CoCalc \cite{LI:SAGE-RZStats}.
Unfortunately samples used for the RV $Xp$ were small ($N=1000$),
due to limitations of using SSAGE in the cloud.
An investigation of the first 100 primes and the pair correlations $c_{pq}=<X_p,X_q>$
revealed this interesting resonance phenomenon,
which deserves a ``proper'' study ($N=1$ million, and a random sample of bigger primes, perhaps).
The present author will focus on a theoretical understanding instead,
in the upcoming work \cite{LI:RHWZ}.

For example, consider the $7$-th prime $p[7]=19$, and plot the correlation coefficient with the
first 100 primes $X_q$ (Fig.5) \cite{LI:PlotXpqv1}.
% X(19) correlations from Plot Xpq correlations v1 save.pdf
\begin{figure}[h!]
\includegraphics[width=6in,height=3in]{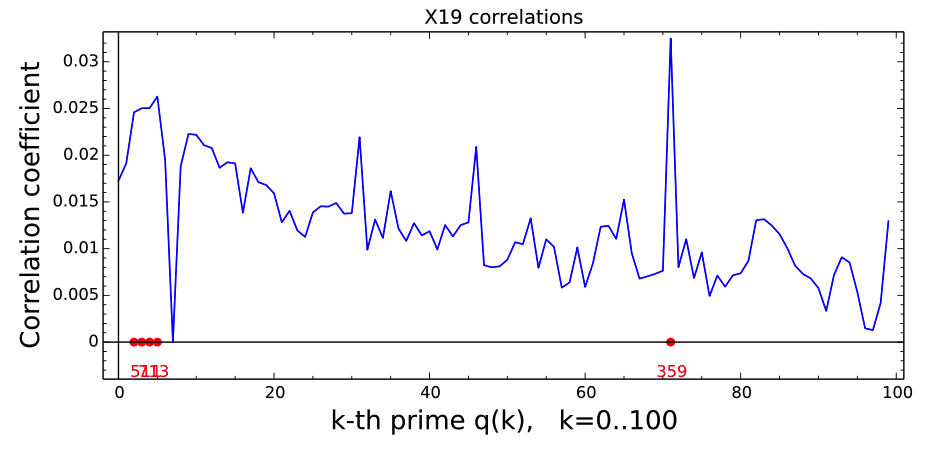}
\caption{$X_{19}$ correlations}
\end{figure}
High correlations occur for small primes $q=3,5,11,13$, for probably trivial reasons (or small sample);
notice for instance the typical envelope of the graph.
High correlation coefficient for big primes occur, in this sample, only  a $q=389$, 
which in its turn, has symmetry factors $q-1=2\cdot 179$ (simple, Saint-Germain prime).

The ``dip'' of the correlation coefficient in the plot occurs at $q=p$, due to resetting the coefficient $c=1$
at zero, for convenience of the plotting (remove the big jump).

As another example, consider $p[9]=29$, and the corresponding correlations of $X_p$ (Fig.6).
\begin{figure}[h!]
\includegraphics[width=6in,height=3in]{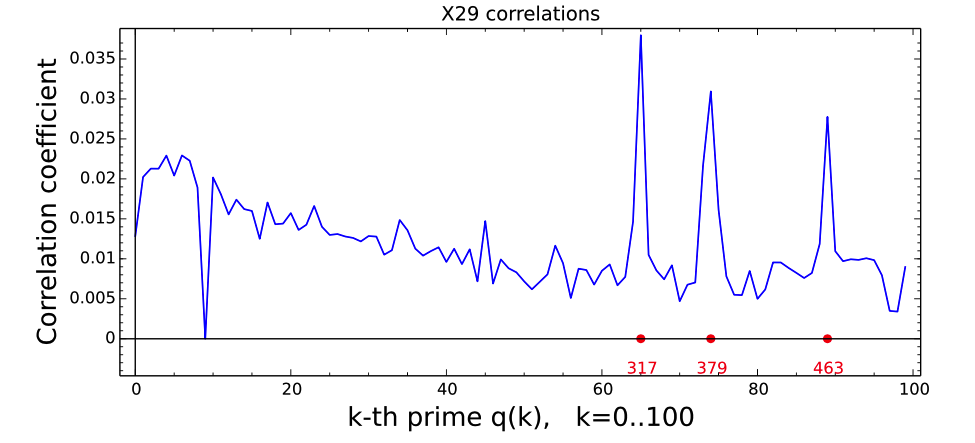}
\caption{$X_{29}$ correlations}
\end{figure}
The ``high resonances'' with $X_q$ occur for $q=p[65]=317$, with symmetries $q-1=2^2\cdot79$,
$p[74]=379$, $q-1=2\cdot3^3\cdot7$, and $q=p[89]=463$, $q-1=2\cdot3\cdot7\cdot11$.

As with the Brahe-Kepler-Newton example of longterm research cycle from {\em data} to {\em math-model},
the theoretical explanation awaits at least a couple of ``laws'' inferred from observation of similar data,
with hindsight from previous attempts \cite{LI:PNRZ,LI:RHWZ}.
To have an idea what theoretical aspect to look for,
we present next some conclusions and a corresponding program of study.

%***************************************************************************************
% 			Conclusions
%***************************************************************************************
\section{Conclusions and Further developments}
% Adeles, Q Algebraic QG, Tannaka-krein {p^n}=Spec Ab_f
After reviewing the evidence for a pro-algebraic structure on the Riemann Spectrum $R-Spec$,
consisting of the imaginary parts of the Riemann zeros,
the statistics by Ford and Zaharescu of the p-sectors and their correlations was recalled \cite{Ford1,Ford2}.

To complement their investigation, the article considers the hierarchy of primes,
represented by the POSet structure \cite{LI:POSet},
and investigates computationally, using SAGE/CoCalc, 
the correlation $c(X_p,X_q)$, between the $p$ and $q$ sectors,
when $q|p-1$, i.e. when multiplication by $q$ is a generator (independent dimension)
of the group of symmetries $Aut_{Ab}(F_p)$ of the ``discrete vector space'' $F_p$
(``elementary finite string'' \cite{LI:Rem}).

This article searches for the $Z-Spec$ POSet shadow on $R-Spec$,
under duality.
A measure of how strong the partial order $q<<p$ on primes 
dualizes to the corresponding random variables (Cramer characters),
is their statistical correlation coefficient $c(X_p,X_q)$.

Several examples show a resonance phenomenon, 
i.e. a sharp significantly higher correlation between $p$ and $q$ sectors of R-Spec
precisely when one prime is a Klein symmetry of the other prime Galois geometry.

Further developments in this direction,
of a computational study, will look to other measures of such a correlation,
for example the Pearson coefficient of correlation or clustering methods.

\vspace{.1in}
Understanding these ``resonances'' awaits for a deeper theoretical explanation,
probably as in physics, in terms of irreducible representations 
of, what else, the absolute Galois group,
at the level of (again, what else!), 
Tannaka-Krein duality  and crystal bases.

A more concrete, preliminary program of study is sketched below.

% The Research Program:
$$\xymatrix{
& \underline{Primes} & \underline{Zeros}\\
Global\dto^{2} & (Z-Spec, POSet) \dto^{Ramification} \rto^{Fourier\ Duality}& (R-Spec,?) \dto{Graal}\\
Local & L_p-functions \rto^{Artin/Weil} &  Weil \ Zeros\\
Categorification & Ab_f \rto^{Tannaka-Krein?} & SL_2(Z)-mod?
}$$
The next step in exploring $R-Spec$ is a review of Weil zeros and their connection
with the Riemann Zeros, beyond a formal analogy \cite{LI:RHWZ}.

%=======================================================================
%                      Bibliography
%=======================================================================

\end{document}